\documentclass{article}

\usepackage{multicol,graphicx,color}
\usepackage{pslatex}
\usepackage{amsthm}
\usepackage{amssymb}

\newcommand{\myfig}[3][0]{
\begin{center}
  \vspace{0.2cm}
  \includegraphics[width=#3\hsize,angle=#1]{#2}

  \nobreak\medskip
\end{center}
}

\newcommand{\mycaption}[1]{
  \vspace{0.2cm}
  \begin{quote}
    {{\sc Figure} \arabic{figure}: #1}
  \end{quote}
  \vspace{0.2cm}
  \stepcounter{figure}
}

\theoremstyle{plain}
\newtheorem{theorem}{Theorem}

\newtheorem*{theorema}{Main Theorem}

\author{Tiancheng Ouyang\footnote{Department of Mathematics, Brigham Young University, Provo, UT 84602, \textit{ouyang@math.byu.edu}}, Skyler C.~Simmons\footnote{Department of Mathematics, Brigham Young University, Provo, UT 84602, \textit{xinkaisen@yahoo.com}}, Duokui Yan\footnote{Chern Institute of Mathematics, Nankai University, Tianjin, 300071, P.R.C., \textit{duokuiyan@gmail.com}}
}

\title{Periodic Solutions with Singularities in Two Dimensions in the $n$-body Problem}

\begin{document}

\maketitle

\begin{abstract}
Analytical methods are used to prove the existence of a periodic,
symmetric solution with singularities in the planar 4-body
problem. A numerical calculation and simulation are used to
generate the orbit. The analytical method easily extends to any
even number of bodies. Multiple simultaneous binary collisions are
a key feature of the orbits generated.
\end{abstract}



\section{Introduction}

The $n$-body problem of celestial mechanics is one of the most important problems in the field of dynamical systems. The following differential equation
\begin{equation}\label{Famous}
m_i \ddot{\rho_i} = \sum_{j\neq i} -\frac{m_i m_j (\rho_i -
\rho_j)}{{\vert \rho_i - \rho_j \vert}^3}
\end{equation}
gives a mathematical description of the planar $n$-body problem, where $\rho_i \in \mathbb{R}^2$ denotes the position of the $i$th body having mass $m_i$. All derivatives are taken with respect to time $t$. The potential energy of the system is given by
\begin{equation}
U = \sum_{1\le i < j \le n} \frac{m_i m_j}{|\rho_i - \rho_j|},
\end{equation}
and the kinetic energy is given by
\begin{equation}
T = \frac{1}{2} \sum_{i=1}^n m_i |\dot{\rho_i}|^2.
\end{equation}\\

Linearly stable symmetric periodic orbits are one aspect of the $n$-body problem. The elliptic Lagrangian triangular periodic orbits are linearly stable for certain values of eccentricity and the three masses \cite{MS},\cite{Ro1}. The Montgomery-Chenciner figure-eight orbit for three equal masses \cite{CM}, \cite{Moore} has been shown by Roberts \cite{GR} to be linearly stable by an innovative symmetry reduction technique he developed.\\

%Here is the info for the two additional references added to the above paragraph

Singularities are another particular aspect of the $n$-body problem. Binary collisions, triple collisions, etc,. are discussed at length in \cite{SM}. The Simultaneous Binary Collision (SBC) problem has been widely studied as well, both analytically and numerically. Sim\'{o} \cite{SI} showed that the block regularization in the cases of the $n$-body problem which reduce to one-dimensional problems is differentiable, but the map passing from initial to final conditions (in suitable choices of transversal sections) is exactly $C^{8/3}$. Ouyang and Yan \cite{OY} give another approach for the regularization and analyze some properties of SBC solutions in the collinear four-body problem. Elbialy \cite{EL} studied the nature of the collision-ejection orbits associated to SBC. \\

Schubart \cite{JS} combined these two aspects to produce a singular linearly stable periodic orbit in the three-body equal mass collinear problem. The motion of the middle mass regularly alternates between binary collisions with each of the outer two masses. His work was subsequently extended to the unequal mass case by both H\'enon \cite{MH} and Hietarinta and Mikkola \cite{HM}. Sweatman \cite{SW} later extended this work to a four-body periodic solution in one dimension, with bodies alternating between SBC of the outer mass pairs and binary collision of the inner two masses.\\

In this paper, we present the analytic existence of a family of singular symmetric periodic planar orbits in the four-body equal mass problem. The initial conditions of these orbits are symmetric in both positions and velocities, which lead to periodic simultaneous binary collisions with each of the four masses alternating between collisions with its two nearest neighbors. Due to the abundance of symmetries present in the initial conditions, we can reduce the number of variables needed to just four -- two for representing position and two for representing momentum. In contrast to its one-dimensional counterparts, the proof for existence of this orbit is surprisingly simple. We begin in Section \textbf{\ref{Describe}} by giving a description of the proposed orbit and prove its existence. In Section \textbf{\ref{Numerical}} we present the numerical methods used to produce the initial conditions that will lead to this orbit. In Section \textbf{\ref{Family}}, we consider variants on this orbit, giving a family of orbits with singularities for an even number of equal masses.\\

Since the initial submission of this paper, we have been doing additional work with Dr. Lennard Bakker (Brigham Young University) and Dr. Gareth Roberts (College of the Holy Cross) implementing Robert's linear stability technique as presented in \cite{GR}. After precisely defining the symmetries that are present in the regularized coordinates, it is shown that the group of symmetries in the orbit is isomorphic to the dihedral group $D_4$. Further, as a consequence of Robert's technique, we have shown that the four-body planar orbit presented in this paper is linearly stable \cite{BORSY}. Further work has also been done on orbits in this family with alternating unequal masses \cite{BOYS}. Rather than a single mass parameter, the bodies have masses $m_1$, $m_2$, $m_1$, $m_2$ as numbered moving counterclockwise through the plane. Since some symmetries have been lost by this change in masses, it is necessary to choose two initial condition parameters as well as two initial velocities. Although numerically this is not a difficult problem, an analytical technique will require much more work.\\

%Here is the info for the addition reference mentioned in the last paragraph

\section{The Proposed Orbit}
\subsection{Analytical Description}\label{Describe}

\begin{center}
\myfig{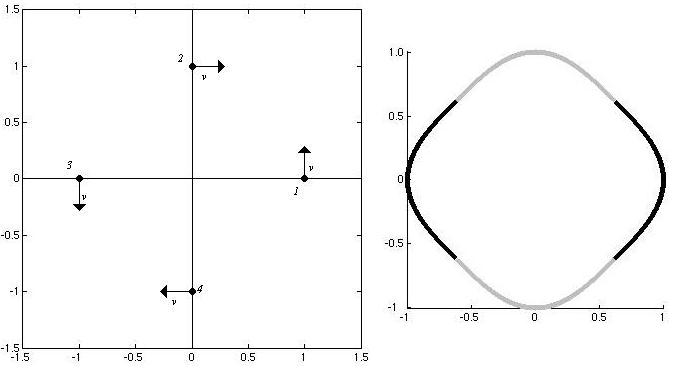}{1} \mycaption{On the left, we illustrate
the initial conditions leading to the four-body two-dimensional
periodic SBC oribt. On the right, the orbit is shown.}
\end{center}

Initially we focused on finding a symmetric, periodic SBC orbit
for four equal masses in two dimensions. Without loss of
generality, we assume that the orbit begins with the four bodies
lying at $(\pm 1,0)$ and $(0, \pm1)$ with initial velocities $(0, \pm v)$ and $(\pm v, 0)$, respectively, where $v \in (0, +\infty)$.  For convenience throughout the rest of the paper, we number the bodies 1 to 4 as in Figure 1.\\

The singularity of SBC in this problem is not essential. For a
better understanding of the behavior of the motion of the bodies
in a neighborhood of a collision, the standard technique is to
make a change of coordinates and rescale time. In the new
coordinates, the orbits which approach collision can be extended
across the collision in a smooth manner with respect to the new
time variable. This technique is called \textit{regularization}.
In our problem, the regularization describes the behavior of the
bodies approaching and escaping collisions,
similar to the collisions of billiard balls.\\

%Regularization is important for studying the collision. Intuitively,
%the singularity at collision can be opened if we consider colliding
%masses with a small but nonzero impact parameter. The bodies will
%rapidly spin around each other. Therefore, if the impact parameter
%approaches zero, we should anticipate that the bodies hit and then
%rebound from the central force in a manner resembling a perfect
%elastic rebound. More significantly, from the mathematical view, if
%collisions can be regularized, collision solutions are
%the limits of nearby non-collision solutions.\\

Due to the symmetry of the initial conditions and the equations
governing the motion of the bodies, the symmetry that is present in
the initial conditions is maintained in the regularized sense.\\

\begin{theorema}
Let $E = T - U$ be the total energy and $m$ be the mass for each
of the four bodies. For any $E<0$ and $m>0$, there exists a
symmetric, periodic, four-body orbit with SBC in $\mathbb{R}^2$.
\end{theorema}

Without loss of generality, we can assume $m=1$ and the initial
positions are as illustrated in Figure 1. The proof will be given at the end of this section.\\

Let $t_0$ be the time of first SBC. For $t \in [0, t_0)$, let the
coordinate of body 1 be $(x_1, x_2)$. By symmetry, the coordinates
of bodies 2, 3, and 4 are $(x_2, x_1)$, $(-x_1, -x_2)$ and $(-x_2,
-x_1)$, respectively. Using equation (\ref{Famous}), the
acceleration of a body at point $(x_1,x_2)$ is given by:
\begin{equation}\label{Newtons}
(\ddot{x_1}, \ddot{x_2}) = - \bigg[
\frac{(x_1-x_2,x_2-x_1)}{(2(x_1-x_2)^2)^\frac{3}{2}} + \frac{(2x_1,
2x_2)}{(4x_1^2+4x_2^2)^\frac{3}{2}} + \frac{(x_1+x_2,
x_1+x_2)}{(2(x_1+x_2)^2)^\frac{3}{2}} \bigg]
\end{equation}

We now perform the regularization of the system. The system has the Hamiltonian:
\begin{equation}
H=\frac{1}{8} (w_1^2+w_2^2)-
\frac{\sqrt{2}}{x_1-x_2}-\frac{\sqrt{2}}{x_1+x_2}-\frac{1}{\sqrt{x_1^2+x_2^2}}
\end{equation}
where $w_1=4 \dot{x}_1$ and $w_2=4 \dot{x}_2$ are the conjugate
momenta to $x_1$ and $x_2$. Note that SBC happens when $x_1= \pm
x_2$. We introduce a new set of coordinates:
$$q_1=x_1-x_2, \qquad q_2= x_1+x_2.$$
Their conjugate momenta $p_i$ are given by using a generating
function $F=(x_1-x_2) p_1 +(x_1+x_2) p_2$:
$$w_1= p_1+p_2, \qquad w_2= p_2-p_1.$$
The Hamiltonian corresponding to the new coordinate system is
\begin{equation}
H= \frac{1}{4}(p_1^2+p_2^2)-
\frac{\sqrt{2}}{q_1}-\frac{\sqrt{2}}{q_2}
-\frac{\sqrt{2}}{\sqrt{q_1^2+q_2^2}}.
\end{equation}

Following the work of Sweatman \cite{SW}, we introduce another
canonical transformation:
$$q_i=Q_i^2, \qquad P_i=2 Q_i p_i \qquad(i=1,2)$$
with $Q_i > 0$.  We also introduce a new time variable $s$, which satisfies
$\frac{d t}{d s}=q_1 q_2$. This produces a regularized Hamiltonian
in extended phase space:
$$\Gamma =\frac{d t}{d s} (H-E)$$
\begin{equation}\label{6}
= \frac{1}{16} (P_1^2Q_2^2 + P_2^2 Q_1^2)- \sqrt{2}(Q_1^2+Q_2^2)
-\frac{\sqrt{2}Q_1^2 Q_2^2 }{\sqrt{Q_1^4+ Q_2^4} } -Q_1^2 Q_2^2 E
\end{equation}
where $E$ is the total energy of the Hamiltonian $H$.\\

The regularized Hamiltonian gives the following differential equations of motion:
\begin{equation}\label{1}
Q_1' = \frac{1}{8}P_1 Q_2^2
\end{equation}
\begin{equation}\label{2}
Q_2' = \frac{1}{8}P_2 Q_1^2
\end{equation}
\begin{equation}\label{3}
P_1' = -\frac{1}{8}P_2^2 Q_1 + 2\sqrt{2}Q_1+\frac{2\sqrt{2}Q_1 Q_2^2}{\sqrt{Q_1^4+Q_2^4}}-\frac{2\sqrt{2}Q_1^5 Q_2^2}{(Q_1^4+Q_2^4)^\frac{3}{2}} + 2EQ_1Q_2^2
\end{equation}
\begin{equation}\label{4}
P_2' = -\frac{1}{8}P_1^2 Q_2 + 2\sqrt{2}Q_2+\frac{2\sqrt{2}Q_2 Q_1^2}{\sqrt{Q_1^4+Q_2^4}}-\frac{2\sqrt{2}Q_2^5 Q_1^2}{(Q_1^4+Q_2^4)^\frac{3}{2}} + 2EQ_2Q_1^2
\end{equation}
with initial conditions
\begin{equation}\label{5}
Q_1(0)=1, \quad Q_2(0)=1,  \quad P_1(0)=-4v,  \quad P_2(0)=4v
\end{equation}
where derivatives are with respect to $s$, and E is the total
energy of the Hamiltonian $H$.\\

\begin{theorem}\label{thm1}
Let $s_0$ be the time of the first SBC in the regularized system.
Then $s_0$ is a continuous function with respect to the initial
velocity $v$. Furthermore,
$$p_2 (t_0)= \frac{P_2 (s_0, v)}{2 Q_2 (s_0, v)}$$
is also continuous with respect to $v$.
\end{theorem}

\begin{proof}
At the first SBC, $Q_1(s_0)= 0$, and $Q_2(s_0) = \sqrt{q_2}=
\sqrt{x_1+x_2} > 0$. Our goal is to show that $p_2 (t_0)$  is a
continuous function with
respect to $v$. \\

Because $\Gamma=0$ at $s=s_0$, $P_1(s_0)= -4 \sqrt[4]{2}$ from
(\ref{6}). Since $\Gamma$ is regularized, the solution $P_i=P_i(s,
v)$ and $Q_i=Q_i(s, v)$ are continuous functions with respect to
the two variables $s$ and $v$. At time $s=s_0$,

$$0= Q_1(s_0(v), v).$$
To apply the implicit function theorem, we need to show that
$$\frac{\partial Q_1}{\partial s} (s_0, v) \neq 0.$$
From (\ref{1})
$$\frac{\partial Q_1}{\partial s} (s_0, v)= \frac{1}{8} P_1 Q_2^2 \mid_{(s_0, v)} =- \frac{1}{2}\sqrt[4]{2} Q_2(s_0) ^2 <0. $$
So $s_0=s_0(v) $ is a continuous function of $v$. Therefore both
$P_2(s_0, v)$ and $Q_2(s_0, v)$ are continuous functions of $v$.
Further, since $Q_2(s_0, v)>0$, $p_2(t_0)$ is also a continuous
function of $v$.
\end{proof}

\begin{theorem}\label{thm2}
There exists a $v=v_0$ such that $\dot{x}_1(t_0)+ \dot{x}_2
(t_0)=\frac{1}{2}p_2 (t_0) = 0$, where $t_0$ is the time of the
first SBC, i.e. the net momentum of bodies 1 and 2 at the first
SBC is 0.
\end{theorem}

The outline of this proof is as follows: We will show that there
exist $v_1$ and $v_2$ such that $\dot{x}_1+\dot{x}_2$ is negative
at SBC for $v=v_1$ and positive at SBC for $v=v_2$. The result
then follows by Theorem \ref{thm1}. \\

%We know that when $v=0$, it is total collision. In this case, $\dot{x}_1+\dot{x}_2$ is $-\infty$ at the time of total collision.\\
\begin{proof}
Consider Newton's equation before the time of the first SBC:

\begin{equation}\label{dx1}
\ddot{x}_1=\frac{x_2-x_1}{2\sqrt{2} (x_1-x_2)^3 } - \frac{2x_1}{8
(x_1^2+x_2^2)^{3/2} } - \frac{x_1+x_2}{2 \sqrt{2} (x_1+x_2)^3 } ,
\end{equation}

\begin{equation}\label{dx2}
 \ddot{x}_2=\frac{x_1-x_2}{2\sqrt{2} (x_1-x_2)^3 } - \frac{2x_2}{8 (x_1^2+x_2^2)^{3/2} } - \frac{x_1+x_2}{2 \sqrt{2} (x_1+x_2)^3 }.
\end{equation}
Therefore,
\begin{equation}\label{y}
\ddot{x}_1+ \ddot{x}_2= - \frac{x_1+x_2}{4 (x_1^2+x_2^2)^{3/2} } -
\frac{1}{ \sqrt{2} (x_1+x_2)^2 } <0,
\end{equation}
which means $\dot{x}_1+ \dot{x}_2$ is decreasing with respect to $t$. \\

At the initial time $t=0$, $x_1=1$, $x_2=0$, $\dot{x}_1=0$, and
$\dot{x}_2= v$. Note that for $v \in (0, \infty)$, there is no
triple collision or total collision for $t \in [0, t_0]$, where
$t_0$ is the time of the first SBC, as a triple collision implies total collapse by symmetry. Also, from the initial time to
$t_0$, $0 \leq x_2 \leq x_1 \leq 1$, $0<x_1+x_2<2$, and
$x_1^2+x_2^2<4$.\\

Let $y(t)= x_1(t)+ x_2(t)$. Then for any choice of $v$,  $
\ddot{y}(t)<0$ and $ 0<y(t)<2$ hold  for any $t \in [0, t_0]$.
In other words, $\dot{y}(t)$ is decreasing with respect to $t$.\\

First, we will show that there exists $v_1$ such that $\dot{y}
(t_0) <0$. When $v=0$ the four bodies form a central configuration
and, as a consequence, the motion of the four bodies leads to
total collapse. Consider the time interval $t \in [0, t_0/2)$. In
this interval, the differential equations (\ref{dx1}) and
(\ref{dx2}) have no singularity, and $\ddot{y}(t_0/2)<0$. By
continuous dependence on initial conditions, $\dot{y}(t_0/2)=
\dot{x}_1 (t_0/2) + \dot{x}_2 (t_0/2) $ is a continuous function
with respect to the initial velocity $v$. When $v=0$, $ \dot{x}_1
(t_0/2) <0 , \ \dot{x}_2 (t_0/2) = 0$, which gives
$\dot{y}(t_0/2)<0 $. Therefore, there exists a $\delta>0$, such
that $\dot{y}(t_0/2)<0$ holds for any $v \in (-\delta,
\delta)$. \\

Choose $v_1= \delta/2$, then $\dot{y}(t_0/2)<0$. Because $
\dot{y}(t)$ is decreasing with respect to $t$, $\dot{y}(t_0) \leq
\dot{y}(t_0/2)<0.$\\

Next we will show that there exists $v_2$ such that  $
\dot{y} (t_0) >0$. Note that as $v \to \infty$,
$$\lim_{v \to \infty} y(t_0)= \lim_{v \to \infty} x_1(t_0)+ x_2(t_0) =2 $$
and
$$\lim_{v \to \infty} \dot{y}(t_0) = \infty.$$
Therefore there exists some positive value $v_2$, such that
$\dot{y}(t_0)>0$.
\end{proof}

\begin{proof}[Proof of the Main Theorem]

From Theorem \ref{thm2}, we know there exists an initial velocity
$v=v_0$ such that $\dot{x}_1(t_0)+\dot{x}_2(t_0)=0$. Let $\{P_1,
P_2, Q_1, Q_2\}$ for $s \in [0, s_0]$ be the solution in the
regularized system corresponding to the orbit from $t=0$ to
$t=t_0$. Following collision, consider the behavior of the first
and second bodies. Assume their velocity was reflected about the
$y=x$ line in the plane. In the new coordinate system, this
corresponds to a new set of functions
$$\{-P_1(2s_0-s), -P_2(2s_0-s), -Q_1(2s_0-s), -Q_2(2s_0-s)\}$$
for $s \in [s_0, 2s_0]$. We can easily check that
$$\{-P_1(2s_0-s), -P_2(2s_0-s), -Q_1(2s_0-s), -Q_2(2s_0-s)\}$$
for $s \in [s_0, 2s_0]$ is also a set of solutions for equations
(\ref{1}) through (\ref{4}) with initial conditions at $s=s_0$.
Also, $\{P_1(s), P_2(s), Q_1(s), Q_2(s)\}$ for $s \in [s_0, 2s_0]$
satisfies equations (\ref{1}) through (\ref{4}) with the same
initial conditions at $s=s_0$. Note that equations (\ref{1})
through (\ref{4}) with initial conditions at $s=s_0$ have a unique
solution for any choice of $v \in (0, \infty)$. Then by
uniqueness, the orbit for $s\in [s_0, 2s_0]$ must be the same as
the orbit for $s \in[0, s_0]$ in reverse, i.e. $$P_i(s)=
-P_i(2s_0-s), Q_i(s)= -Q_i(2s_0 -s)$$ for $s \in[0, s_0]$.
Therefore at time $s=2s_0$, bodies 1 and 2 will have returned to
their initial positions with velocities $(0,-v)$ and $(-v, 0)$
respectively. Similarly, at time $s=2s_0$, bodies 3 and 4 will
have also returned to their initial positions with velocities $(0,
v)$ and $(v,0)$ respectively.\\

Next, we use symmetry and uniqueness to show the orbit from
$s=2s_0$ to $s=4s_0$ and the orbit from $s=0$ to $s=2s_0$ will be
symmetric with respect to the $y-$axis. Compare the motion of body
2 and body 3 from $s=2s_0$ to $s=4s_0$ with the motion of body 2
and body 1 from time $s=0$ to $s=2s_0$. The initial conditions of
body 3 at $s=2s_0$ and the initial conditions of body 1 at $s=0$
are symmetric with respect to the $y-$axis. Also the initial
conditions of body 2 at $s=2s_0$ and the initial conditions of
body 4 at $s=0$ are symmetric with respect to the $x-$axis.
Therefore, by uniqueness, the orbit of bodies 2 and 3 from
$s=2s_0$ to $s=4s_0$ and the orbit of bodies 1 and 2 from $s=0$ to
$s=2s_0$ must be symmetric with respect to $y-$axis. Therefore,
the orbit of bodies 1 and 4 from $s=2s_0$ to $s=4s_0$ and and the
orbit of bodies 3 and 4 from $s= 0$ to $s=2s_0$ are symmetric with
respect to the $y-$axis.  Hence, at $s=4s_0$, the positions and velocities of the four
bodies are exactly the same as at $s=0$. Therefore, the orbit  is
periodic with period $s=4s_0$.
\end{proof}

It is worth noting here that the previous proof implies a
time-reversing symmetry for the periodic orbit. This provides
further evidence for the conjecture made by Roberts \cite{GR},
stating that linearly stable periodic orbits in the equal mass
$n$-body problem must have a time-reversing symmetry.  (Linear
stability of this orbit is shown in \cite{BORSY}.)

\subsection{Numerical Method}\label{Numerical}
As we are searching for a periodic orbit of the $n$-body problem,
we assume the value of the Hamiltonian needs to be negative. Using the
initial positions of the four bodies described earlier, it is not
hard to find the potential energy at $t=0$: $$U = {2\sqrt{2} +
1}.$$ Then, acting under the negative Hamiltonian assumption:
$${2\sqrt{2}+1} \geq \sum_{i=1}^n \frac{m_i |v_i|^2}{2}.$$ Since
all masses are equal, if we require that the velocities of each
body are equal in magnitude, we obtain:
\begin{equation}
v_{max} = \sqrt{\frac{2\sqrt{2}+1}{2}}
\end{equation}
with $v_{max}$ defined to be the value of $v$ such that the value
of the Hamiltonian is zero. Define $\theta = \frac{v}{v_{max}}$.
This parameter is used in the numerical algorithm. \\

At this point it becomes necessary to find out just how much
kinetic energy is required to obtain the periodic orbit. Since we
know suitable bounds on the velocity parameter ($\theta \in (0,
1)$), we can search the interval numerically. We use an $n$-body
simulator with the initial positions previously described. The
simulation is ran until one SBC occurs. For simplicity, we
consider only the collision between the first and second bodies in
the first quadrant. Summing their velocities immediately before
the collision gives a vector running along the line $y=x$ (due to
symmetry), with both components having the same sign.  The
magnitude of this vector is given in Figure 2.  Negative
magnitudes represent vectors with both components less than zero.

\begin{center}
\myfig{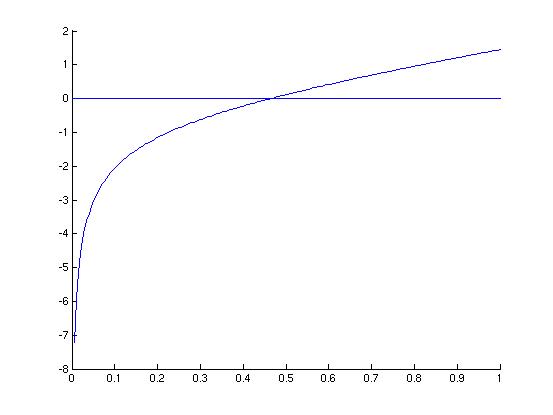}{.75} \mycaption{The magnitude of the net
velocity of the first two bodies (vertical axis) at the time of
collision for various values of $\theta$ (horizontal axis)}.
\end{center}

Next, a standard bisection method is used to find the amount of
energy required to cause the net velocity at collision to be zero.
Using the initial interval $\theta \in [0, 1]$ and iterating to a
tolerance of $10^{-14}$, the correct value of $\theta$ was found
to be $\theta = 0.46449539554694$. \\

It is worth noting that both the proof of existence and the numerical method do not guarantee the uniqueness of this orbit.  Numerical simulations demonstrate that for values of $\theta$ near the correct value, the orbit remains for a significant length of time with the paths of the bodies lying in a ``fattened'' annular region roughly the shape of the original orbit.  Near the extreme ends, the orbit experiences near total-collapse and fall apart rapidly.  Although we do not focus on these questions just yet, a more thorough study of the dynamics could prove to be quite interesting.

\section{Variants}\label{Family}

%\subsection{Orbits of arbitrary radius}

%The orbit described in the previous section is not the only orbit in
%this family. If we let the four bodies have initial position $\pm R$
%on the $x$ and $y$ axes, the potential energy at $t=0$ is given by:
%$$U=\frac{2\sqrt{2}+1}{R}$$
%This yields the maximum velocity for each body:
%\begin{equation}
%v_{max} = \sqrt{\frac{2\sqrt{2}+1}{2R}}
%\end{equation}

%Again, for each individual value of $R$, we let $\theta =
%\frac{v}{v_{max}}$. The following figure shows the magnitudes of
%the net velocity between the first and second bodies at collision
%time for various values of $R$.

%\begin{center}
%\myfig{RVVN04.jpg}{.75} \mycaption{Curves showing the magnitude of
%the net velocity of the first two bodies (vertical axis) at the
%time of collision for various values of $\theta$ (horizontal axis)
%for $R=0.5, 1.0, 1.5, 2.0, 2.5$ near the zero point. Again,
%negative magnitudes represent vectors with both components less
%than zero. The zeros of the curves do not occur at the same value
%of $\theta$ on each curve, but the numerical difference is minor
%($<0.01$).}
%\end{center}

\subsection{Orbits of more than four bodies}

The same technique can be adopted to find similar orbits for any
arbitrary even number $n$. A key feature of these orbits will be
higher numbers of simultaneous binary collisions. For a given
value of $n$, initial positions are given by spacing the bodies
evenly about the unit circle. The potential energy (and the value
of $v_{max}$) is found numerically by iterating over each pair of
planets and summing the reciprocal of the distances between them.
(Recall that all $m_i = 1$.)  Velocities are then assigned to the
bodies in alternating counter-clockwise and clockwise directions,
initially tangent to the circle. Again we consider the collision
between the first and the second bodies. Although the net velocity
of the two at collision will not lie along the $y=x$ line, the
components of this vector will both have the same sign. The
magnitudes of the net velocity between the first two bodies at
initial collision are shown in Figure 3 for various values of $n$.
Lower curves in the graph correspond to higher values of $n$.
Again, negative magnitudes correspond to both components being
negative. \\

\begin{center}
\myfig{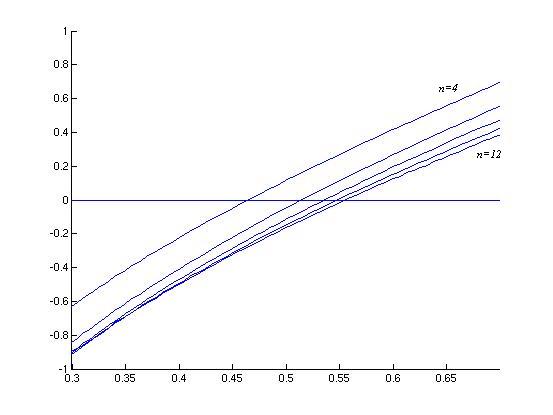}{.75} \mycaption{Curves showing the magnitude of
the net velocity of the first two bodies (vertical axis) at the
time of collision for various values of $\theta$ (horizontal axis)
for $n=4, 6, 8, 10, 12$.}
\end{center}

Pictures of the orbit for $n=6$ and $n=8$ are shown in Figure 4.
It is readily seen that as $n$ increases, the shape of the orbit
more closely approximates a circle.

\begin{center}
\myfig{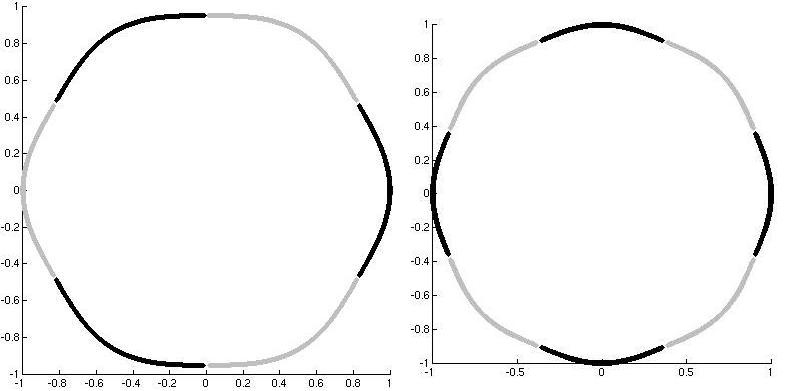}{.75} \mycaption{The six- and eight-body
two-dimensional periodic SBC orbits.}
\end{center}

\noindent{\textbf{Acknowledgements}}

Skyler Simmons would like to give thanks to Dr. Kening Lu, who
recommended me for this project and special thanks to Jason Grout
for hours of discussion that helped point my mathematical career
in the right direction. We are also all indebted to the advice
given by Dr. Lennard Bakker and all the referees during the
revision process.\\

%MAKE BOLD TEXT LIKE \textbf{THIS}, and italicized text like \it{THIS}.]

\end{document}